\documentclass[12pt]{amsart}

\usepackage{amsmath}

\allowdisplaybreaks[4]

\newtheorem{theorem}{\bf Theorem}

\newtheorem{lemma}{\bf Lemma}
\newtheorem{proposition}{\bf Proposition}

\theoremstyle{definition}
\newtheorem{definition}{\bf Definition}

\theoremstyle{remark}

\newcommand{\Z}{\ensuremath{\mathbb{Z}}}

\newcommand{\R}{\ensuremath{\mathbb{R}}}

\newcommand{\N}{\ensuremath{\mathbb{N}}}

\newcommand{\sumln}{{\mbox{$\sum\limits_{l=1\,}^N$}}}
\newcommand{\ip}[2]{{\mbox{$\left<#1,#2\right>$}}}

\begin{document}

\title[Wavelet packets on local fields] {Wavelet packets and wavelet frame packets on local fields}

\author{Biswaranjan Behera}

\address{(B. Behera) Theoretical Statistics and Mathematics Unit, Indian Statistical Institute, 203 B. T. Road, Kolkata, 700108, India}

\email{biswa@isical.ac.in}

\author{Qaiser Jahan}
\address{(Q. Jahan) Theoretical Statistics and Mathematics Unit, Indian Statistical Institute, 203 B. T. Road, Kolkata, 700108, India}
\email{qaiser\_r@isical.ac.in}

\thanks{Research of the second author is supported by a grant from CSIR, India.}



\date{\today}




\begin{abstract}
Using a prime element of a local field $K$ of positive characteristic $p$, the concepts of multiresolution analysis (MRA) and wavelet can be generalized to such a field. We prove a version of the splitting lemma for this setup and using this lemma we have constructed the wavelet packets associated with such MRAs. We show that these wavelet packets generate an orthonormal basis by translations only. We also prove an analogue of splitting lemma for frames and construct the wavelet frame packets in this setting.
\end{abstract}

\maketitle

\section{Introduction}
The concepts of wavelet and multiresolution analysis of $\R^n$ has been extended to many different setups. Dahlke~\cite{Dah} introduced it in locally compact groups (see also~\cite{Far, La1, La2, La3}). It was generalized to abstract Hilbert spaces by Han, Larson, Papadakis, Stavropoulos~\cite{HLPS, SP}. Lemarie~\cite{Le} extended this concept to stratified Lie groups.

Recently, R. L. Benedetto and J. J. Benedetto~(\cite{BB, Ben}) developed a wavelet theory for local fields and related groups. Albeverio, Kozyrev, Khrennikov, Shelkovich, Skopina and their collaborators also discuss about MRA and wavelets on the $p$-adic field $\mathbb{Q}_p$ in a series of papers~\cite{AK, KS, KSS, Koz}. Note that $\mathbb{Q}_p$ is a local field of characteristic $0$.

Jiang, Li and Jin~\cite {JLJ} gave a definition of MRA on a local field of positive characteristic $p$ and, similar to $\R^n$, have constructed the wavelets from an MRA.

In this article we construct the wavelet packets associated with such an MRA. We also generalize the wavelet frame packets to this setup. First of all, we will discuss about wavelet packets very briefly.

Let $\{V_j:j\in\Z\}$ be an MRA of $L^2(\R)$ with scaling function $\varphi$ and wavelet $\psi$. Let $W_j$ be the corresponding wavelet subspaces: $W_j=\overline{\rm span}\{\psi_{jk}:k\in\Z\}$. In the construction of a wavelet from an MRA, essentially the space $V_1$ is split into two orthogonal components $V_0$ and $W_0$. Note that $V_1$ is the closure of the linear span of the  functions $\{2^{1/2}\varphi(2\cdot-k):k\in\Z\}$, whereas $V_0$ and $W_0$  are respectively the closure of the span of $\{\varphi(\cdot-k):k\in\Z\}$ and $\{\psi(\cdot-k):k\in\Z\}$. Since $\varphi(2\cdot-k)=\varphi\left(2(\cdot-2^{-1}k)\right)$, we see that the above procedure splits the half integer translates of a function into integer translates of two functions.

In a similar way, we can split $W_j$, which is the span of  $\{\psi(2^j\cdot-k):k\in\Z\}= \{\psi\left(2^j(\cdot-2^{-j}k)\right):k\in\Z\}$, to get two functions whose $2^{-(j-1)}k$ translates will span the same space $W_j$. Repeating the splitting procedure $j$ times, we get $2^j$ functions whose integer translates alone  span the space $W_j$. If we apply this to each $W_j$, then the resulting basis of $L^2(\R)$ will consist of integer translates of a countable number of functions (instead of all dilations and translations of the wavelet $\psi$). This basis is called the ``wavelet packet basis". The concept of wavelet packets was introduced by Coifman, Meyer and Wickerhauser~\cite{cmw1,cmw2}. For a nice exposition of wavelet packets of $L^2(\R)$ with dilation 2, we refer to~\cite{hw}.

The concept of wavelet packet was subsequently generalized to $\R^n$ by taking tensor products~\cite{cm}. The non-tensor product versions are due to Shen~\cite{she} for dyadic dilation, and Behera~\cite{Beh} for MRAs with a general dilation matrix and several scaling functions. Other notable generalizations are the biorthogonal wavelet packets~\cite{cd}, non-orthogonal version of wavelet packets~\cite{cl}, the wavelet frame packets~\cite{che} on $\R$ for dilation 2, and the orthogonal, biorthogonal and frame wavelet packets on $\R^n$ by Long and Chen~\cite{lc} for the dyadic dilation.

We have organized the article as follows. In section~2, we discuss some preliminary facts about local fields. In section~3, we introduce the concept of MRA on a local field $K$ of positive characteristic and prove a crucial lemma called the splitting lemma. We construct the wavelet packets in section~4 and prove that they generate an orthonormal basis for $L^2(K)$. In section~5, we prove some basic results needed to prove an analogue of splitting lemma for wavelet frames on $K$ and finally in section~6 the wavelet frame packets are constucted.

\section{Preliminaries on local fields}

Let $K$ be a field and a topological space. Then $K$ is called a \emph{locally compact field} or a \emph{local field} if both $K^+$ and $K^*$ are locally compact abelian groups, where $K^+$ and $K^*$ denote the additive and multiplicative groups of $K$ respectively.

If $K$ is any field and is endowed with the discrete topology, then $K$ is a local field. Further, if $K$ is connected, then $K$ is either $\R$ or $\mathbb{C}$. If $K$ is not connected, then it is totally disconnected. So by a local field, we mean a field $K$ which is locally compact, nondiscrete and totally disconnected.

We use the notation of the book by Taibleson~\cite{Taib}. Proofs of all the results stated in this section can be found in the books~\cite{Taib} or~\cite{RV}.

Let $K$ be a local field. Since $K^+$ is a locally compact abelian group, we choose a Haar measure $dx$ for $K^+$. If $\alpha\neq 0, \alpha\in K$, then $d(\alpha x)$ is also a Haar measure. Let $d(\alpha x) = |\alpha|dx$. We call $|\alpha|$ the \emph{absolute value} or \emph{valuation} of $\alpha$. We also let $|0| = 0$.

The map $x\rightarrow |x|$ has the following properties:
\begin{itemize}
\item[(a)] $|x| = 0$ if and only if $x=0$;
\item[(b)] $|xy|=|x||y|$ for all $x,y\in K$;
\item[(c)] $|x+y|\leq\max\{|x|, |y| \}$ for all $x,y\in K$.
\end{itemize}
Property (c) is called the \emph{ultrametric inequality}.

The set $\mathfrak{D}=\{x\in K : |x|\leq 1\}$ is called the \emph{ring of integers} in $K$. It is the unique maximal compact subring of $K$. Define $\mathfrak{P}=\{x\in K:|x|<1\}$. The set $\mathfrak{P}$ is called the \emph{prime ideal} in $K$. The prime ideal in $K$ is the unique maximal ideal in $\mathfrak{D}$. It is principal and prime.

Since $K$ is totally disconnected, the set of values $|x|$ as $x$ varies over $K$ is a discrete set of the form $\{s^k: k\in\Z\}\cup \{0\}$ for some $s>0$. Hence, there is an element of $\mathfrak{P}$ of maximal absolute value. Let $\mathfrak{p}$ be a fixed element of maximum absolute value in $\mathfrak{P}$. Such an element is called a \emph{prime element} of $K$. Note that as an ideal in $\mathfrak{D},\mathfrak{P}=\left\langle \mathfrak{p}\right\rangle = \mathfrak{p}\mathfrak{D}$.

It can be proved that $\mathfrak{D}$ is compact and open. Hence, $\mathfrak{P}$ is compact and open. Therefore, the residue space $\mathfrak{D}/\mathfrak{P}$ is isomorphic to a finite field $GF(q)$, where $q=p^c$ for some prime $p$ and $c\in\N$. For a proof of this fact we refer to ~\cite{Taib}.

For a measurable subset $E$ of $K$, let $|E|=\int_K\chi_E(x)dx$, where $\chi_E$ is the characteristic function of $E$ and $dx$ is the Haar measure of $K$ normalized so that $|\mathfrak{D}|=1$. Then, it is easy to see that $|\mathfrak{P}|=q^{-1}$ and $|\mathfrak{p}|=q^{-1}$ (see~\cite{Taib}). It follows that if $x\neq 0$, and $x\in K$, then $|x|=q^k$ for some $k\in\mathbb{Z}$.

Let $\mathfrak{D}^*=\mathfrak{D}\setminus\mathfrak{P}=\{x\in K: |x|=1\}$. $\mathfrak{D}^*$ is the group of units in $K^*$. If $x\neq 0$, we can write $x=\mathfrak{p}^k x'$, with $x'\in\mathfrak{D}^*$.

Recall that $\mathfrak{D}/\mathfrak{P}\cong GF(q)$. Let $\mathcal{U}=\{a_i\}_{i=0}^{q-1}$ be any fixed full set of coset representatives of $\mathfrak{P}$ in $\mathfrak{D}$. Let $\mathfrak{P}^k=\mathfrak{p}^k\mathfrak{D}=\{x\in K: |x|\leq q^{-k}\}, k\in\Z$. These are called \emph{fractional ideals}. Each $\mathfrak{P}^k$ is compact and open and is a subgroup of $K^+$ (see~\cite{RV}). Then, if $x\in\mathfrak{P}^k, k\in\Z$, $x$ can be expressed uniquely as $x=\sum_{l=k}^{\infty} c_l\mathfrak{p}^l, c_l\in \mathcal{U}$.

If $K$ is a local field, then there is a nontrivial, unitary, continuous character $\chi$ on $K^+$. It can be proved that $K^+$ is self dual (see~\cite{Taib}).

Let $\chi$ be a fixed character on $K^+$ that is trivial on $\mathfrak{D}$ but is nontrivial on $\mathfrak{P}^{-1}$. We can find such a character by starting with any nontrivial character and rescaling. We will define such a character for a local field of positive characteristic. For $y\in K$, we define $\chi_y(x)=\chi(yx)$, $x\in K$.

\begin{definition}
If $f\in L^1(K)$, then the Fourier transform of $f$ is the function $\hat f$ defined by
\[
\hat f(\xi)= \int_K f(x)\overline{\chi_{\xi}(x)}~dx.
\]
\end{definition}
Note that
\[
\hat f(\xi)= \int_K f(x)\overline{\chi(\xi x)}~dx = \int_K f(x)\chi(-\xi x)~dx.
\]

Similar to the standard Fourier analysis on the real line, one can prove the following results.
\begin{itemize}
\item[(a)] The map $f\rightarrow \hat f$ is a bounded linear transformation of $L^1(K)$ into $L^{\infty}(K)$, and $\|\hat f\|_{\infty}\leq \|f\|_1$.
\item[(b)] If $f\in L^1(K)$, then $\hat{f}$ is uniformly continuous.
\item[(c)] $f\in L^1 \cap L^2(K)$, then $\|\hat {f} \|_2=\|f\|_2$.
\end{itemize}

To define the Fourier transform of function in $L^2(K)$, we introduce the functions $\Phi_k$. For $k\in\Z$, let $\Phi_k$ be the characteristic function of $\mathfrak{P}^k$.
\begin{definition}
For $f\in L^2(K)$, let $f_k=f\Phi_{-k}$ and
\[
\hat{f}(\xi)=\lim\limits_{k\rightarrow\infty}\hat f_k(\xi)
=\lim\limits_{k\rightarrow \infty}\int_{\left|x\right|\leq q^k} f(x)\overline{\chi_{\xi}(x)}~d\xi,
\]
where the limit is taken in $L^2(K)$.
\end{definition}
We have the following theorem (see Theorem 2.3 in~\cite{Taib}).
\begin{theorem}
The fourier transform is unitary on $L^2(K)$.
\end{theorem}

Let $\chi_u$ be any character on $K^+$. Since $\mathfrak{D}$ is a subgroup of $K^+$, the restriction $\chi_{u}|_\mathfrak{D}$ is a character on $\mathfrak{D}$. Also, as character on $\mathfrak{D}, \chi_u = \chi_v$ if and only if $u-v\in \mathfrak{D}$. That is, $\chi_u=\chi_v$ if $u+\mathfrak{D}=v+\mathfrak{D}$ and $\chi_u\neq \chi_v$ if $(u+\mathfrak{D})\cap (v+\mathfrak{D})=\phi$. Hence, if $\{u(n)\}_{n=0}^{\infty}$ is a complete list of distinct coset representative of $\mathfrak{D}$ in $K^+$, then $\{\chi_{u(n)}\}_{n=0}^{\infty}$ is a list of distinct characters on $\mathfrak{D}$. It is proved in~\cite{Taib} that this list is complete. That is, we have the following proposition.
\begin{proposition}
Let $\{u(n)\}_{n=0}^{\infty}$ be a complete list of (distinct) coset representatives of $\mathfrak{D}$ in $K^+$. Then $\{\chi_{u(n)}\}_{n=0}^{\infty}$ is a complete list of (distinct) characters on $\mathfrak{D}$. Moreover, it is a complete orthonormal system on $\mathfrak{D}$.
\end{proposition}
Given such a list of characters $\{\chi_{u(n)}\}_{n=0}^{\infty}$, we define the Fourier coefficients of $f\in L^1(\mathfrak{D})$ as
\[
\hat{f}(u(n))=\int_{\mathfrak{D}}f(x)\overline{\chi_{u(n)}(x)}dx.
\]
The series $\sum\limits_{n=0}^{\infty}\hat{f}(u(n))\chi_{u(n)}(x)$ is called the Fourier series of $f$. From the standard $L^2$ theory for compact abelian groups we conclude that the Fourier series of $f$ converges to $f$ in $L^2(\mathfrak{D})$ and
\[
\int_{\mathfrak{D}}|f(x)|^2dx= \sum\limits_{n=0}^{\infty}|\hat{f}(u(n))|^2.
\]
Also, if $f\in L^1(\mathfrak{D})$ and $\hat f(u(n))=0$ for all $n\in\N_0$, then $f=0$ a. e.

These results hold irrespective of the ordering of the characters. We now proceed to impose a natural order on the sequence $\{u(n)\}_{n=0}^{\infty}$. Note that $\Gamma=\mathfrak{D}/\mathfrak{P}$ is isomorphic to the finite field $GF(q)$ and $GF(q)$ is a $c$-dimensional vector space over the field $GF(p)$. We choose a set $\{1=\epsilon_0, \epsilon_1, \epsilon_2, \cdots, \epsilon_{c-1}\} \subset\mathfrak{D}^*$ such that span$\{\epsilon_j\}_{j=0}^{c-1}\cong GF(q)$.
Let $\N_0=\N\cup \{0\}$. For $n\in \N_0$ such that $0\leq n< q$, we have
\[
n=a_0+a_1 p+\cdots+a_{c-1} p^{c-1},\quad 0\leq a_k<p, k=0,1,\cdots,c-1.
\]
Define
\begin{equation}\label{e.undef}
u(n)=(a_0+a_1\epsilon_1+\cdots+a_{c-1}\epsilon_{c-1})\mathfrak{p}^{-1}.
\end{equation}
Now, write
\[
n=b_0+b_1q+b_2q^2+\cdots+b_sq^s,\quad 0\leq b_k<q, k=0,1,2,\cdots,s,
\]
and define
\[
u(n)=u(b_0)+u(b_1)\mathfrak{p}^{-1}+\cdots+u(b_s)\mathfrak{p}^{-s}.
\]
Note that $u(0)=0$ and $\{u(n)\}_{n=0}^{q-1}$ is a complete set of coset representatives of $\mathfrak{D}$ in $\mathfrak{P}^{-1}$ (see~\cite{Taib}). Hence, $\{u(n)\mathfrak{p}\}_{n=0}^{q-1}$ is a complete set of coset representatives of $\mathfrak{P}$ in $\mathfrak{D}$. Therefore,
\[
\{u(n)\mathfrak{p}\}_{n=0}^{q-1}\cong \mathfrak{D}/\mathfrak{P}\cong GF(q)\cong {\rm span}\{\epsilon_j\}_{j=0}^{c-1}.
\]

In general, it is not true that $u(m+n)=u(m)+u(n)$. But
\begin{equation}\label{eq.un}
u(rq^k+s)=u(r)\mathfrak{p}^{-k}+u(s)\quad{\rm if}~r\geq 0, k\geq 0~{\rm and}~0\leq s <q^k.
\end{equation}
For brevity, we will write $\chi_n=\chi_{u(n)}, n\geq 0$. As mentioned before, $\{\chi_n\}_{n=0}^{\infty}$ is a complete set of characters on $\mathfrak{D}$.

Let $\mathcal{U}=\{a_i\}_{i=0}^{q-1}$ be a fixed set of coset representatives of $\mathfrak{P}$ in $\mathfrak{D}$. Then every $x\in K$ can be expressed uniquely as
\[
x=x_0+\sum\limits_{k=1}^nb_k\mathfrak{p}^{-k}, \quad x_0\in \mathfrak{D}, b_k\in \mathcal{U}.
\]

Let $K$ be a local field characteristic $p>0$ and $\epsilon_0, \epsilon_1, \dots, \epsilon_{c-1}$ be as above. We define a character $\chi$ on $K$ as follows:
\begin{equation}\label{chi}
\chi(\epsilon_{\mu}\mathfrak{p}^{-j})=
\left\{
\begin{array}{lll}
\exp(2\pi i/p), & \mu=0~\mbox{and}~j=1,\\
1, & \mu=1,\cdots,c-1~\mbox{or}~j\neq 1.
\end{array}
\right.
\end{equation}
Note that $\chi$ is trivial on $\mathfrak{D}$ but nontrivial on $\mathfrak{P}^{-1}$.

In order to be able to define the concepts of multiresolution analysis and wavelets on local fields, we need analogous notions of translation and dilation. Since $\bigcup\limits_{j\in\Z}\mathfrak{p}^{-j}\mathfrak{D}=K,$ we can regard $\mathfrak{p}^{-1}$ as the dilation (note that $|\mathfrak{p}^{-1}|=q$) and since $\{u(n): n\in\N_0\}$ is a complete list of distinct coset representatives of $\mathfrak{D}$ in $K$, the set $\{u(n): n\in\N_0\}$ can be treated as the translation set. So we make the following definition.
\begin{definition}
A finite set $\{\psi_m:m=1,2,\cdots,M\}\subset L^2(K)$ is called a \emph{set of basic wavelets} of $L^2(K)$ if the system $\{q^{j/2}\psi_m(\mathfrak{p}^{-j}\cdot-u(k)): 1\leq m\leq M, j\in\Z, k\in\N_0\}$ forms an orthonormal basis for $L^2(K)$.
\end{definition}


\section{Multiresolution analysis on local fields and the splitting lemma}

Similar to $\R^n$, wavelets can be constructed from a multiresolution analysis which we define below (see~\cite{JLJ}).
\begin{definition}
Let $K$ be a local field of characteristic $p>0$, $\mathfrak{p}$ be a prime element of $K$ and $u(n)\in K$ for $n\in \mathbb{N}_0$ be as defined above. A multiresolution analysis (MRA) of $L^2(K)$ is a sequence $\{V_j\}_{j\in\mathbb{Z}}$ of closed subspaces of $L^2(K)$ satisfying the following properties:
\begin{enumerate}
\item[(a)] $V_j\subset V_{j+1}$ for all $j\in\Z$;
\item[(b)] $\bigcup\limits_{j\in\Z}V_j$ is dense in $L^2(K)$ and $\bigcap\limits_{j\in\Z}V_j = \{0\}$;
\item[(c)] $f\in V_j$ if and only if $f(\mathfrak{p}^{-1}\cdot)\in V_{j+1}$ for all $j\in\Z$;
\item[(d)] there is a function $\varphi\in V_0$, called the \emph{scaling function}, such that $\{\varphi(\cdot-u(k)): k\in\N_0\}$ forms an orthonormal basis for $V_0$.
\end{enumerate}
\end{definition}

Given an MRA $\{V_j: j\in\Z\}$, we define another sequence $\{W_j: j\in\Z\}$ of closed subspaces of $L^2(K)$ by
\[
W_j=V_{j+1}\ominus V_j.
\]
These subspaces also satisfy
\begin{equation}\label{e.wj}
f\in W_j~ {\rm if~and~only~if} ~f(\mathfrak{p}^{-1}\cdot)\in W_{j+1},~ j\in\Z.
\end{equation}
Moreover, they are mutually orthogonal, and we have the following orthogonal decompositions:
\begin{eqnarray}
L^2(K) & = & \bigoplus\limits_{j\in\Z}W_j \label{e.L2decomp1} \\
& = & V_0 \oplus\Bigl(\bigoplus\limits_{j\geq 0}W_j\Bigr).
\end{eqnarray}

Observe that the dilation is induced by $\mathfrak{p}^{-1}$ and $\left|\mathfrak{p}^{-1}\right|=q$. As in the case of $\R^n$, we expect the existence of $q-1$ number of functions $\{\psi_1,\psi_2,\cdots,\psi_{q-1}\}$ to form a set of basic wavelets. In view of ~(\ref{e.wj}) and (\ref{e.L2decomp1}), it is clear that if $\{\psi_1,\cdots,\psi_{q-1}\}$ is a set of function such that $\{\psi_m(\cdot-u(k)): 1\leq m\leq M, k\in\N_0\}$ forms an orthonormal basis for $W_0$, then $\{q^{j/2}\psi_m(\mathfrak{p}^{-j}\cdot-u(k)): 1\leq m\leq M, j\in\Z, k\in\N_0\}$ forms an orthonormal basis for $L^2(K)$.

For $f\in L^2(K)$, we define
\[
f_{j,k}=q^{j/2}f(\mathfrak{p}^{-j}x-u(k)), \quad j\in\Z, k\in\N_0.
\]
Then it is easy to see that
\[
\|f_{j,k}\|_2=\|f\|_2
\]
and
\[
(f_{j,k})^\wedge (\xi)= q^{-j/2}\overline{\chi_k(\mathfrak{p}^j\xi)}\hat{f}(\mathfrak{p}^j\xi).
\]

Since $\varphi\in V_0\subset V_1$, and $\{\varphi_{1,k}:k\in\N_0\}$ is an orthonormal basis in $V_1$, we have
\[
\varphi(x)=\sum\limits_{k\in\N_0}h_kq^{1/2}\varphi(\mathfrak{p}^{-1}x-u(k)),
\]
where $h_k=\langle\varphi,\varphi_{1,k}\rangle$ and $\{h_k:k\in\N_0\}\in\ell^2(\N_0)$. Taking Fourier transform, we get
\begin{eqnarray}
\hat{\varphi}(\xi) & = & q^{-1/2}\sum\limits_{k\in\N_0}h_k\overline{\chi_k(\mathfrak{p}\xi)}\hat{\varphi}(\mathfrak{p}\xi) \nonumber \\
& = & m_0(\mathfrak{p}\xi)\hat{\varphi}(\mathfrak{p}\xi) \label{e.m0},
\end{eqnarray}
where $m_0=q^{-1/2}\sum\limits_{k\in\N_0}h_k\overline{\chi_k(\xi)}$.

Let us call a function $f$ on $K$ \emph{integral-periodic} if
\[
f(x+u(k))=f(x)~\mbox{for all}~k\in\N_0.
\]
The following facts were proved in~\cite{JLJ}.
\begin{itemize}
\item[(a)] $\chi_k(u(l))=\chi(u(k)u(l))=1$ for all $k, l\in\N_0$.
\item[(b)] The function $m_0$ is integral-periodic.
\item[(c)] The system $\{\varphi(\cdot-u(k)): k\in\N_0\}$ is orthonormal if and only if $\sum\limits_{k\in\N_0}\left|\widehat{\varphi}(\xi+u(k))\right|^2=1$ a.e.
\end{itemize}

Given an MRA of $L^2(K)$, suppose that there exist $q-1$ integral-periodic functions $m_l$, $1\leq l\leq q-1$, such that the matrix
\[
M(\xi)=\Big(m_l(\mathfrak{p}\xi+\mathfrak{p}u(k))\Big)_{l,k=0}^{q-1}
\]
is unitary. It was also proved in~\cite{JLJ} that $\{\psi_1,\psi_2,\cdots,\psi_{q-1}\}$ is a set of basic wavelets of $L^2(K)$ if we define
\[
\hat{\psi}_l(\xi)=m_l(\mathfrak{p}\xi)\hat{\varphi}(\mathfrak{p}\xi).
\]

We now prove a lemma, the splitting lemma, which is essential for the construction of wavelet packets. With the help of this lemma, we can decompose a closed subspace of $L^2(K)$ into finitely many mutually orthogonal subspaces in a suitable manner.

\begin{lemma}[The splitting lemma]\label{lem:split}
Let $\varphi \in L^2(K)$ be such that $\{\varphi(\cdot-u(k)): k\in \mathbb{N}_0\}$ is an orthonormal system.
Let $V = \overline{\rm span}\{q^{1/2}\varphi(\mathfrak{p}^{-1}\cdot-u(k)): k\in \mathbb{N}_0\}$. Let
$m_l = q^{-1/2}\sum _{k=0}^{\infty}h_{k}^{l} \overline{\chi_k}(\xi)$, $0\leq l\leq q-1$, where
$\{h_{k}^{l}: k\in \mathbb{N}_0\}\in\ell^2(\mathbb{N}_0)$ for $0\leq l\leq q-1$. Define $\hat\psi _{l}(\xi) = m_l(\mathfrak{p}\xi)\hat{\varphi}(\mathfrak{p}\xi)$. Then $\{\psi_l (\cdot - u(k)): 0\leq l\leq q-1, k\in \N_0\}$ is an orthonormal system in $V$ if and only if the matrix
\begin{equation*}
M(\xi) = \Bigl(m_l(\mathfrak{p}\xi + \mathfrak{p}u(k))\Bigr)_{l,k=0}^{q-1}
\end{equation*}
is unitary for a.e $\xi\in\mathfrak{D}$.

Moreover, $\{\psi_l(\cdot - u(k)): 0\leq l\leq q-1 , k\in\N_0\}$ is an orthonormal basis of $V$ whenever it is orthonormal.
\end{lemma}

\proof Assume that $M(\xi)$ is unitary for a.e $\xi\in\mathfrak{D}$. Then, for $0\leq s, t\leq q-1$ and $k, l\in\N_0$, we have
\begin{eqnarray*}
\lefteqn{\Bigl\langle\psi_s\bigl(\cdot - u(k)\bigr),\psi_t\bigl(\cdot - u(l)\bigr)\Bigr\rangle} \\
& = & \Bigl\langle\Bigl(\psi_s\bigl(\cdot - u(k)\bigr)\Bigr)^\wedge,
\Bigl(\psi_t\bigl(\cdot - u(l)\bigr)\Bigr)^\wedge\Bigr\rangle \\
& = & \int_K\hat\psi_s(\xi)\overline{\chi_k(\xi)}\overline{\hat{\psi}_t(\xi)}\chi_l(\xi)~d\xi \\
& = & \int_{\mathfrak{D}}\sum\limits_{n\in\N_0}\hat\psi_s\bigl(\xi+u(n)\bigr)
    \overline{\hat\psi_t\bigl(\xi+u(n)\bigr)}\overline{\chi_k(\xi)}\chi_l(\xi)~d\xi \\
& = & \int_{\mathfrak{D}}\sum\limits_{n\in\N_0}m_s\bigl(\mathfrak{p}\xi+\mathfrak{p}u(n)\bigr)
    \overline{m_t\bigl(\mathfrak{p}\xi+\mathfrak{p}u(n)\bigr)} \bigl|\hat\varphi\bigl(\mathfrak{p}\xi+\mathfrak{p}u(n)\bigr)\bigr|^2\overline{\chi_k(\xi)}\chi_l(\xi)~d\xi \\
& = & \int_{\mathfrak{D}}\sum\limits_{\mu=0}^{q-1}\sum\limits_{n\in\N_0}m_s
    \bigl(\mathfrak{p}\xi+\mathfrak{p}u(qn+\mu)\bigr)
    \overline{m_t\bigl(\mathfrak{p}\xi+\mathfrak{p}u(qn+\mu)\bigr)} \\
&   & \qquad \times\bigl|\hat\varphi\bigl(\mathfrak{p}\xi+\mathfrak{p}u(qn+\mu)\bigr)\bigr|^2
    \overline{\chi_k(\xi)}\chi_l(\xi)~d\xi \\
& = & \int_{\mathfrak{D}}\Big\{\sum\limits_{\mu=0}^{q-1}m_s(\mathfrak{p}\xi+\mathfrak{p}u(\mu))
    \overline{m_t(\mathfrak{p}\xi+\mathfrak{p}u(\mu))}\Big\}\overline{\chi_k(\xi)}\chi_l(\xi)~d\xi \\
& = & \int_{\mathfrak{D}}\delta_{s,t}\overline{\chi_k(\xi)}\chi_l(\xi)~d\xi \\
& = & \delta_{s,t}\delta_{k,l}.
\end{eqnarray*}
Hence, $\{\psi_s(\cdot-u(k)): 0\leq s \leq q-1, k\in\N_0\}$ is an orthonormal system in $V$. The converse can be proved by reversing the above steps.

To prove the second part, let $f\in V$ be such that $f$ is orthogonal to $\psi_l(\cdot-u(k))$ for all $l=0, 1,\dots,q-1$, $k\in\N_0$. We claim that $f=0$ a.\ e.

Since $f \in V$, we have
\[
f(x) = \sum\limits_{m\in\N_0}q^{1/2}c_m\varphi(\mathfrak{p}^{-1}x-u(m)),
\]
for some $\{c_m:m\in\N_0\}\in\ell^2(\N_0)$. So there exists an integral periodic function  $m_f$ in $L^2(\mathfrak{D})$ such that
\[
\hat{f}(\xi) = m_f(\mathfrak{p}\xi)\hat{\varphi}(\mathfrak{p}\xi).
\]

Hence, for all $l=0, 1,\dots,q-1$, $k\in\N_0$, we have (by a similar calculation)
\begin{eqnarray*}
0 & = & \bigl\langle f,\psi_l(\cdot-u(k))\bigr\rangle\\
& = & \int_K\hat{f}(\xi)\overline{\hat{\psi}_l(\xi)}\chi_k(\xi)d\xi \\
& = & \int_{\mathfrak{D}}\Bigl\{\sum\limits_{\mu=0}^{q-1}m_f(\mathfrak{p}\xi+\mathfrak{p}u(\mu))
\overline{m_l(\mathfrak{p}\xi+\mathfrak{p}u(\mu))}\Bigr\}\chi_k(\xi)d\xi.
\end{eqnarray*}
Therefore, for all $l=0, 1,\dots,q-1$, we have
\begin{equation*}
\sum\limits_{\mu=0}^{q-1}m_f\bigl(\mathfrak{p}\xi + \mathfrak{p}u(\mu)\bigr) \overline{m_l\bigl(\mathfrak{p}\xi+\mathfrak{p}u(\mu)\bigr)} = 0.
\end{equation*}
Now, for a.e. $\xi$, the vector $\Bigl(m_f(\mathfrak{p}\xi+\mathfrak{p}u(\mu))\Bigr)_{\mu=0}^{q-1} \in \mathbb{C}^q$, being orthogonal to each row vector of the unitary matrix $M(\xi)$, is the zero vector. In particular, $m_f(\mathfrak{p}\xi)=0$ a.e. This means $\hat f=0$ a. e. and hence $f=0$ a. e.
\qed


\section{Construction of wavelet packets}

Let $\{V_j:j\in\Z\}$ be an MRA of $L^2(K)$ and $\varphi$ be the corresponding scaling function. Then we have (see~\eqref{e.m0}),
\[
\hat{\varphi}(\xi)=m_0(\mathfrak{p}\xi)\hat{\varphi}(\mathfrak{p}\xi).
\]
Applying the splitting lemma to $V=V_1$, we get $\{\psi_l(\cdot-u(k)): 0\leq l \leq q-1,k\in\N_0\}$ is an orthonormal basis for $V_1$. Now we will define a sequence $\{\omega_n: n\geq 0\}$ of functions. Let
\[
\omega_0 = \varphi
\]
and
\[
\omega_n = \psi_n \quad(1\leq n \leq q-1),
\]
where
\begin{equation}\label{e.hatpsil}
\hat{\psi_l}(\xi) = m_l(\mathfrak{p}\xi)\hat{\varphi}(\mathfrak{p}\xi) \quad(1\leq l \leq q-1).
\end{equation}
Suppose $\omega_m$ is defined for $m\geq 0$. For $0 \leq r \leq q-1$, define
\begin{equation}\label{wpkt}
\omega_{r+qm}(x) = q^{1/2}\sum\limits_{k\in \N_0}h_{k}^r \omega_m(\mathfrak{p}^{-1}x-u(k)).
\end{equation}
Note that this defines $\omega_n$ for every integer $n\geq 0$. Taking Fourier Transform, we have
\begin{equation}\label{e.ftwpkt}
(\omega_{r+qm})^\wedge (\xi) = m_r(\mathfrak{p}\xi)\hat{\omega}_m(\mathfrak{p}\xi).
\end{equation}

\begin{definition}
The functions $\{\omega_n: n\geq 0,\}$ as defined above will be called the \emph{wavelet packets} corresponding to the MRA $\{V_j: j\in\Z\}$ of $L^2(K)$.
\end{definition}

In the following proposition we find an expression for the Fourier transforms of the wavelet packets in terms of $\hat\varphi$.
\begin{proposition}
For an integer $n\geq 1$, consider the unique expansion of $n$ in the base $q$:
\begin{equation}\label{eqn:q.ary}
n=\mu_1+\mu_2q+\mu_3q^2+\cdots+\mu_jq^{j-1},
\end{equation}
where $0\leq\mu_i\leq q-1$ for all $i=1,2,\dots ,j$ and $\mu_j\not=0$.
Then
\begin{equation}
\hat{\omega}_n(\xi)= m_{\mu_1}(\mathfrak{p}\xi)m_{\mu_2}(\mathfrak{p}^2\xi)\cdots m_{\mu_j}(\mathfrak{p}^j\xi)\hat{\varphi}(\mathfrak{p}^j\xi).
\end{equation}
\end{proposition}
\proof
We will prove it by induction. If $n$ has an expansion as in \eqref{eqn:q.ary}, then we say that it has length $j$.
Since $\omega_0=\varphi$, and $\omega_n=\psi_n, 1\leq n \leq q-1$, it follows from~\eqref{e.m0} and \eqref{e.hatpsil} that the claim is true for length 1. Assume that it is true for length $j$. Let $m$ be an integer with an expansion of length $j+1$. Then there exist integers $\gamma_1, \gamma_2, \dots, \gamma_{j+1}$ with $0\leq\gamma_1, \gamma_2, \dots, \gamma_{j+1}\leq q-1$ and $\gamma_{j+1}\not=0$ such that
\begin{eqnarray*}
m & = & \gamma_1+\gamma_2q+\cdots+\gamma_jq^{j-1}+\gamma_{j+1}q^j \\
  & = & \gamma_1+kq,
\end{eqnarray*}
where $k=\gamma_2+\gamma_3q+\cdots+\gamma_{j+1}q^{j-1}$. Note that $k$ has length $j$. Hence,
\begin{eqnarray*}
\hat{\omega}_m(\xi) & = & (\omega_{\gamma_1+kq})^\wedge(\xi)  \\
& = & m_{\gamma_1}(\mathfrak{p}\xi)\hat{\omega}_k(\mathfrak{p}\xi) \qquad \qquad~({\rm by~\eqref{e.ftwpkt}})  \\
& = & m_{\gamma_1}(\mathfrak{p}\xi)m_{\gamma_2}(\mathfrak{p}^2\xi)\cdots m_{\gamma_j+1}(\mathfrak{p}^{j+1}\xi)\hat{\varphi}(\mathfrak{p}^{j+1}\xi).
\end{eqnarray*}
Hence the induction is complete.
\qed

We will prove the following theorem regarding the wavelet packets.
\begin{theorem}\label{thm:wavpack}
Let $\{\omega_n:n\geq 0\}$ be the basic wavelet packets constructed above. Then
\begin{enumerate}
\item [(i)]
$\{\omega_n(\cdot-u(k)): q^j\leq n \leq q^{j+1}-1,k\in\N_0\}$ is an orthonormal basis of $W_j$, $j\geq 0$.
\item [(ii)]
$\{\omega_n(\cdot -u(k)): 0\leq n \leq q^j-1,k\in\N_0\}$ is an orthonormal basis of $V_j$, $j\geq 0$.
\item [(iii)]
$\{\omega_n(\cdot -u(k)): n \geq 0,k\in\N_0\}$ is an orthonormal basis of $L^2(K)$.
\end{enumerate}
\end{theorem}
\proof Since $\{\omega_n:1\leq n\leq q-1\}$ are wavelets, the case $j=0$ in (i) is trivial. Assume (i) for $j$. We will prove for $j+1$. By our assumption, the set of functions $\{\omega_n(\cdot -u(k)): q^j\leq n \leq q^{j+1}-1,k\in\N_0\}$ is an orthonormal basis of $W_j$. By \eqref{e.wj}, we have
\[
\{q^{1/2}\omega_n(\mathfrak{p}^{-1}\cdot -u(k)): q^j\leq n \leq q^{j+1}-1,k\in\N_0\}
\]
is an orthonormal basis for $W_{j+1}$. Let
\[
E_n = \overline{\rm span}\{q^{1/2}\omega_n(\mathfrak{p}^{-1}\cdot -u(k)): k\in\N_0\}.
\]
Hence,
\begin{equation}\label{eqn:wjp1}
W_{j+1}=\bigoplus\limits_{n=q^j}^{q^{j+1}-1} E_n.
\end{equation}
Applying the splitting lemma to $E_n$, we get the functions
\begin{equation*}
g^{n}_l(x) = \sum\limits_{k=0}^{\infty} h^l_{k}q^{1/2}\omega_n(\mathfrak{p}^{-1}x-u(k)),\quad (0\leq l\leq q-1)
\end{equation*}
such that $\{g^{n}_l(\cdot-u(k)):0\leq l\leq q-1,k\in\N_0\}$ forms an orthonormal basis for $E_n$. Hence,
$\{g^{n}_l(\cdot-u(k)):0\leq l\leq q-1,q^j \leq n \leq q^{j+1}-1,k\in\N_0\}$ forms an orthonormal basis for $W_{j+1}$.
But by~(\ref{wpkt}),
\[
g^{n}_l = \omega_{l+qn}.
\]
This fact, together with (\ref{eqn:wjp1}), shows that
\begin{eqnarray*}
\lefteqn{\{\omega_{l+qn}(\cdot-u(k)):0\leq l\leq q-1,~ q^j\leq n \leq q^{j+1}-1,~k\in\N_0\} }\\
& = & \{\omega_n(\cdot -u(k)): q^{j+1}\leq n \leq q^{j+2}-1,~ k\in\N_0\}
\end{eqnarray*}
is an orthonormal basis for $W_{j+1}$. So (i) is proved. Item (ii) follows from the observation that $V_j = V_0\oplus W_0\oplus\cdots\oplus W_{j-1}$ and (iii) follows from the fact that $\overline{\cup V_j} = L^2(K)$.
\qed

\section{Wavelet frame packets}
Let $\mathcal H$ be a separable Hilbert space. A sequence $\{x_k: k\in\Z\}$ of $\mathcal H$ is said to be a frame for $\mathcal H$ if there exist constants $C_1$ and $C_2$, $0<C_1\leq C_2<\infty$ such that for all $x\in \mathcal H$
\begin{equation}\label{eqn:frame}
C_1\|x\|^2 \leq\sum\limits_{k\in\Z} |\left<x,x_k\right>|^2 \leq C_2\|x\|^2.
\end{equation}
The largest $C_1$ and the smallest $C_2$ for which (\ref{eqn:frame}) holds are called the frame bounds.

Suppose that $\Phi=\{\varphi_1, \varphi_2,\dots, \varphi_N\}\subset L^2(K)$ be such that the system $\{\varphi_l(\cdot-u(k)): 1\leq l\leq N,k\in\N_0\}$ is a frame for its closed linear span $S(\Phi)$. Let $\psi_1,\psi_2,\dots,\psi_N$ be elements in $S(\Phi)$. A natural question to ask is: When can we say that $\left\{\psi_l(\cdot-u(k)):1\leq l\leq N, k\in\N_0 \right\}$ is also a frame for $S(\Phi)$?

If $\psi_j\in S(\Phi)$, then there exists $\left\{p_{jlk}: k\in\N_0\right\}$ in $\ell^2(\N_0)$ such that
\[
\psi_j(x)=\sumln\sum\limits_{k\in\N_0} p_{jlk}\varphi_l(x-u(k)).
\]
Taking Fourier transform, we get
\begin{eqnarray*}
\hat\psi_j(\xi) & = & \sumln\sum\limits_{k\in\N_0} p_{jlk}\overline{\chi_k(\xi)}\hat\varphi_l(\xi) \\
                & = & \sumln P_{jl}(\xi)\hat\varphi_l(\xi),
\end{eqnarray*}
where $P_{jl}(\xi)=\sum\limits_{k\in\N_0} p_{jlk}\overline{\chi_k(\xi)}$. Let $P(\xi)$ be the $N\times N$ matrix:
\[
P(\xi)=\Bigl(P_{jl}(\xi)\Bigr)_{1\leq j,l\leq N}.
\]

Let $S$ and $T$ be two positive definite matrices of order $N$. We say $S\leq T$ if $T-S$ is positive definite. The following lemma is the generalization of Lemma 3.1 in~\cite{che}.
\begin{lemma}\label{lem:frame}
Let $\varphi_l, \psi_l $ for $1\leq l\leq N$, and $P(\xi)$ be as above. Suppose that there exist constants $C_1$ and $C_2$, $0< C_1\leq C_2 < \infty$ such that
\begin{equation}\label{eqn:pstar}
C_1I\leq P^*(\xi)P(\xi) \leq C_2I \quad for~a.e.~\xi\in \mathfrak{D}.
\end{equation}
Then, for all $f\in L^2(K)$, we have
\begin{eqnarray}\label{eqn:phipsi}
\lefteqn{C_1 \sum\limits_{l=1}^N\sum\limits_{k\in\N_0} \left|\ip{f}{\psi_l(\cdot-u(k))}\right|^2\leq \sumln \sum\limits_{k\in\N_0}\left| \ip{f}{\varphi_l(\cdot-u(k))}\right|^2} \\
&  & \leq C_2 \sumln \sum\limits_{k\in\N_0} \left|\ip{f}{\varphi_l(\cdot-u(k))}\right|^2.\nonumber
\end{eqnarray}
\end{lemma}

\proof
For $f, g\in L^2(K)$, we define
\[
\left[f,g\right](\xi) = \sum\limits_{l\in\N_0}\widehat{f}(\xi+u(l))\overline{\widehat{g}(\xi+u(l))}.
\]
Then, for $f\in L^2(K)$, we have
\begin{eqnarray*}
\left[f,\psi_j\right](\xi) & = & \sum\limits_{l\in\N_0}\hat{f}(\xi+u(l))\overline{\hat{\psi}_l(\xi+u(l))} \\
 & = & \sum\limits_{k=1}^{N}\sum\limits_{l\in\N_0}\overline{P_{jk}(\xi+u(l))}\hat{f}(\xi+u(l))
 \overline{\hat{\varphi}_k(\xi+u(l))} \\
 & = & \sum\limits_{k=1}^N \overline{P}_{jk}(\xi)\left[f,\varphi_k\right](\xi),
\end{eqnarray*}
since $P_{jk}$ are integral periodic function. Hence
\[
\sum\limits_{j=1}^N\left|\left[f,\psi_j\right]\right|^2 = \sum\limits_{k,k'=1}^N\sum\limits_{j=1}^N \overline{P}_{jk}P_{jk'}\left[f,\varphi_k\right]\overline{\left[f,\varphi_{k'}\right]}
= XP^*PX^*,
\]
where
\[
X = \bigl(\left[f,\varphi_1\right],\cdots\left[f,\varphi_n\right]\bigr).
\]
By Plancherel Theorem,
\[
\sum\limits_{k\in\N_0}\sum\limits_{l=1}^N\left|\left\langle f,\varphi_l(\cdot-u(k))\right\rangle\right|^2 = \sum\limits_{l=1}^N\int_{\mathfrak{D}}\left|\left[f,\varphi_l\right](\xi)\right|^2d\xi.
\]
Hence, inequality~\eqref{eqn:phipsi} is equivalent to
\[
C_1\int_{\mathfrak{D}}XX^* \leq \int_{\mathfrak{D}}XP^*PX^* \leq C_2\int_{\mathfrak{D}}XX^*,
\quad{\rm for~all}~f\in L^2(K).
\]
This follows from~\eqref{eqn:pstar}.
\qed

We now introduce a matrix $E(\xi)$. For $0\leq r,s\leq q-1$ and $1\leq l,j\leq N$, define for a.e. $\xi$
\begin{equation*}
{\mathcal E}^{rs}_{lj}(\xi)=\delta_{lj}q^{-\frac{1}{2}}\overline{\chi\bigl(u(r)(\xi+\mathfrak{p}u(s))\bigr)}.
\end{equation*}
Let
\begin{equation*}
E^{rs}(\xi)=\Bigl({\mathcal E}^{rs}_{lj}(\xi)\Bigr)_{1\leq l,j\leq N}
\end{equation*}
and
\begin{equation}\label{eqn:matE}
E(\xi)=\Bigl(E^{rs}(\xi)\Bigr)_{0\leq r,s\leq q-1}.
\end{equation}
So $E(\xi)$ is a block matrix with $q$ blocks in each row and each column, and each block is a square matrix of order $N$, so that $E(\xi)$ is a square matrix of order $qN$.

We have the following lemma which will be useful for the splitting trick for frames. In the first part of the lemma we use a technique used by Zheng in~\cite{Z}.
\begin{lemma}\label{lem:matrix}
\begin{enumerate}
\item[(i)] For $0\leq r,s\leq q-1$,
\[
{\frac{1}{q}}\sum\limits_{t=0}^{q-1}\chi\bigl((u(r)-u(s))\mathfrak{p}u(t)\bigr) = \delta_{r,s}.
\]
\item[(ii)] The matrix $E(\xi)$, defined in {\rm(\ref{eqn:matE})}, is unitary for a.e. $\xi\in\mathfrak{D}$.
\end{enumerate}
\end{lemma}
\proof
(i) If $r=s$ then $u(r)-u(s)=0$, hence the left hand side equals 1. We assume $r\neq s$. Let
\[
r=a_0+a_1p+\cdots+a_{c-1}p^{c-1}~{\rm and}~s=b_0+b_1p+\cdots+b_{c-1}p^{c-1}
\]
where $0\leq a_j, b_j\leq p-1$ for $j=0, 1, \dots, c-1$. Then (see~\eqref{e.undef})
\[
u(r)\mathfrak{p}=a_0\epsilon_0+a_1\epsilon_1+\cdots+a_{c-1}\epsilon_{c-1}~{\rm and}~ u(s)\mathfrak{p}=b_0\epsilon_0+b_1\epsilon_1+\cdots+b_{c-1}\epsilon_{c-1}.
\]
Now, let
\[
t=d_0+d_1p+\cdots+d_{c-1}p_{c-1}, 0\leq d_j\leq p-1~\mbox{for}~j=0, 1, \dots, c-1.
\]
Observe that as $t$ varies from $0$ to $q-1$, the integers $d_0, d_1, \dots, d_{c-1}$ all vary from $0$ to $p-1$. For each $j=0, 1, \dots, c-1$, we write
\[
u(r)\mathfrak{p}\epsilon_j=
\gamma_{r,0}^j\epsilon_0+\gamma_{r,1}^j\epsilon_1+\cdots+\gamma_{r,c-1}^j\epsilon_{c-1}
\]
for some unique $\gamma_{r,l}^j\in GF(p), 0\leq l\leq c-1$. Similarly,
\[
u(s)\mathfrak{p}\epsilon_j=
\gamma_{s,0}^j\epsilon_0+\gamma_{s,1}^j\epsilon_1+\cdots+\gamma_{s,c-1}^j\epsilon_{c-1}
\]
for some unique $\gamma_{s,l}^j\in GF(p), 0\leq l\leq c-1$. By the definition of the character $\chi$ (see~\eqref{chi}), we have
\[
\chi(u(r)\mathfrak{p}u(t))=\exp\big(\tfrac{2\pi i}{p}(\gamma_{r,0}^0d_0+\cdots+\gamma_{r,0}^{c-1}d_{c-1})\big)
\]
and
\[
\chi(u(s)\mathfrak{p}u(t))=\exp\big(\tfrac{2\pi i}{p}(\gamma_{s,0}^0d_0+\cdots+\gamma_{s,0}^{c-1}d_{c-1})\big).
\]
Therefore,
\begin{eqnarray*}
\lefteqn{
\sum\limits_{t=0}^{q-1}\chi\bigl((u(r)-u(s))\mathfrak{p}u(t)\bigr) }\\
& = & \sum\limits_{t=0}^{q-1}\chi\bigl(u(r)\mathfrak{p}u(t)\bigr)\overline{\chi\big(u(s)\mathfrak{p}u(t)\bigr)} \\
& = & \sum\limits_{d_0=0}^{p-1}\cdots \sum\limits_{d_{c-1}=0}^{p-1}
      \exp\Big(\tfrac{2\pi i}{p} (\gamma_{r,0}^0d_0+\cdots+\gamma_{r,0}^{c-1}d_{c-1})\Big) \\
&   & \qquad \qquad
      \exp\Big(\tfrac{-2\pi i}{p}(\gamma_{s,0}^0d_0+\cdots+\gamma_{s,0}^{c-1}d_{c-1})\Big) \\
& = & \Bigg(\sum\limits_{d_0=0}^{p-1}\exp\Bigl(\tfrac{2\pi i}{p}(\gamma_{r,0}^0-\gamma_{s,0}^0)d_0\Bigr)\Bigg)\cdots
      \Bigg(\sum\limits_{d_{c-1}=0}^{p-1}\exp\Bigl(\tfrac{2\pi i}{p}(\gamma_{r,0}^{c-1}-\gamma_{s,0}^{c-1})d_{c-1}\Bigr)\Bigg).
\end{eqnarray*}

Since $r\ne s$, we claim that $\gamma_{r,0}^j\not=\gamma_{s,0}^j$ for some $j$, $0\leq j\leq c-1$. If $\gamma_{r,0}^j=\gamma_{s,0}^j$ for all $j$,  then, since $u(r)\mathfrak{p}\neq u(s)\mathfrak{p}$, we have
\begin{eqnarray*}
GF(q) & = & {\rm span}\{(u(r)\mathfrak{p}-u(s)\mathfrak{p})\epsilon_j\}_{j=0}^{c-1}\\
      & = & {\rm span}\big\{(\gamma_{r,0}^j-\gamma_{s,0}^j)\epsilon_0, \cdots, (\gamma_{r,c-1}^j-\gamma_{s,c-1}^j) \epsilon_{c-1}\big\}_{j=0}^{c-1}\\
      & = & {\rm span}\{\epsilon_1,\epsilon_2,\cdots,\epsilon_{c-1}\}.
\end{eqnarray*}
This is a contradiction which proves the claim. Now for any $j$ such that $\gamma_{r,0}^j\not=\gamma_{s,0}^j$, we have
\[
\sum\limits_{d_j=0}^{p-1}\exp\Bigl(\tfrac{2\pi i}{p}(\gamma_{r,0}^j-\gamma_{s,0}^j)d_j\Bigr)
=\tfrac{1-\exp\bigl(2\pi i(\gamma_{r,0}^j-\gamma_{s,0}^j)\bigr)}{1-\exp\bigl(\tfrac{2\pi i}{p}(\gamma_{r,0}^j-\gamma_{s,0}^j)\bigr)}=0,
\]
since $\gamma_{r,0}^j-\gamma_{s,0}^j$ is an integer with $|\gamma_{r,0}^j-\gamma_{s,0}^j|<p$. This proves (i).

To prove (ii), observe that the $(r,s)$-th block of the matrix $E(\xi)E^*(\xi)$ is
\[
\sum\limits_{t=0}^{q-1} E^{rt}(\xi)\left(E^{ts}(\xi)\right)^*.
\]
The $(l,j)$-th entry in this block is
\begin{eqnarray*}
& = & \sum\limits_{t=0}^{q-1}\sum\limits_{m=0}^{N} {\mathcal E}^{rt}_{lm}(\xi)
        \left({\mathcal E}^{ts}_{mj}(\xi)\right)^* \\
& = &  \sum\limits_{t=0}^{q-1}\sum\limits_{m=0}^{N}\delta_{lm} q^{-1/2}\overline{\chi\bigl(u(r)(\xi+\mathfrak{p}u(t))\bigr)}
       \cdot\delta_{jm} q^{-1/2}\chi\bigl(u(s)(\xi+\mathfrak{p}u(t))\bigr) \\
& = & \sum\limits_{m=1}^N\delta_{lm}\delta_{jm}q^{-1}\sum\limits_{t=0}^{q-1}\overline{\chi\bigl(u(r)
       (\xi+\mathfrak{p}u(t))\bigr)}\chi\bigl(u(s)(\xi+\mathfrak{p}u(t))\bigr) \\
& = & \sum\limits_{m=1}^N\delta_{lm}\delta_{jm}\chi((u(s)-u(r))\xi)
        q^{-1}\sum\limits_{t=0}^{q-1}\chi\bigl((u(s)-u(r))\mathfrak{p}u(t))\bigr) \\
& = & \sum\limits_{m=1}^N\delta_{lm}\delta_{jm}\delta_{rs},\quad{\rm (by~part~(i)~of~the~lemma)} \\
& = & \delta_{lj}\delta_{rs}.
\end{eqnarray*}
Hence $E(\xi)E^*(\xi)=I$. Similarly, $E(\xi)^*E(\xi)=I$. Therefore, $E(\xi)$ is a unitary matrix.
\qed

\section{Splitting lemma for wavelet frame packets}

Let $\{\varphi_j: 1\leq j\leq N\}$ be functions in $L^2(K)$ such that $\{\varphi_j(\cdot-u(k)): 1\leq j \leq N,k\in\N_0\}$ is a frame for its closed linear span $V$. For $1\leq l \leq N$, $0\leq r\leq q-1$, suppose that there exist sequences $\{h^r_{ljk}:k\in\Z\}\in\ell^2(\mathfrak{D})$. Define
\[
\psi^r_l(x) = q^{1/2}\sum\limits_{j=1}^N \sum\limits_{k\in\N_0} h^r_{ljk}\varphi_j(\mathfrak{p}^{-1}x-u(k)).
\]
Taking Fourier transform, we get
\[
\hat\psi_l^r(\xi)=\sum\limits_{j=1}^N \sum\limits_{k\in\N_0} h^r_{ljk}q^{-1/2}\overline{\chi_k(\mathfrak{p}\xi)}\hat\varphi_j(\mathfrak{p}\xi)
=\sum\limits_{j=1}^N h^r_{lj}\hat\varphi_j(\mathfrak{p}\xi),
\]
where,
\[
h^r_{lj}(\xi)=\sum\limits_{k\in\N_0}q^{-1/2}h^r_{ljk}\overline{\chi_k(\xi)}.
\]
Let
\[
H_r(\xi)=\bigl(h^r_{lj}(\xi)\bigr)_{1\leq l,j\leq N}
\]
and
\[
H(\xi)=\Bigl(H_r(\xi+\mathfrak{p}u(s))\Bigr)_{0\leq r,s \leq q-1}.
\]
Note that $H(\xi)$ is a square matrix of order $qN$. We can write $\psi^r_l$ as
\begin{eqnarray*}
\psi^r_l(x) & = & \sum\limits_{j=1}^N \sum\limits_{k\in\N_0} h^r_{ljk}q^{1/2}\varphi_j\bigl(\mathfrak{p}^{-1}x-u(k)\bigr) \\
& = & \sum\limits_{j=1}^N \sum\limits_{s=0}^{q-1} \sum\limits_{k\in\N_0}h^r_{lj,qk+s}q^{1/2}\varphi_j\bigl(\mathfrak{p}^{-1}x-u(qk+s)\bigr) \\
& = & \sum\limits_{j=1}^N \sum\limits_{s=0}^{q-1} \sum\limits_{k\in\N_0}h^r_{lj,qk+s}\varphi_j^{(s)}(x-u(k)),
\end{eqnarray*}
where
\begin{equation} \label{eqn:phisj}
\varphi_j^{(s)}(x)= q^{1/2}\varphi_j(\mathfrak{p}^{-1}x-u(s)), \quad0\leq s \leq q-1.
\end{equation}
Note that $u(qk+s)=\mathfrak{p}^{-1}u(k)+u(s)$ (see eq.~\eqref{eq.un}). Taking Fourier transform, we obtain
\begin{eqnarray*}
(\psi^r_l)^\wedge(\xi) & = & \sum\limits_{j=1}^N \sum\limits_{s=0}^{q-1} \sum\limits_{k\in\N_0}h^r_{lj,qk+s}\overline{\chi_{k}(\xi)}(\varphi^{(s)}_j)^\wedge (\xi) \\
& = & \sum\limits_{j=1}^N \sum\limits_{s=0}^{q-1} p^{rs}_{lj}(\xi)(\varphi^{(s)}_j)^\wedge (\xi), \\
\end{eqnarray*}
where $p^{rs}_{lj}(\xi)=\sum\limits_{k\in\N_0} h^r_{lj,qk+s}\overline{\chi_k(\xi)}$. Define the matrices
\begin{equation*}
P^{rs}(\xi)=\Bigl(p^{rs}_{lj}(\xi)\Bigr)_{1\leq l,j\leq N}.
\end{equation*}
and
\begin{equation*}
P(\xi)=\Bigl(P^{rs}(\xi)\Bigr)_{0\leq r,s\leq q-1}.
\end{equation*}
\begin{proposition}
$H(\xi)=P(\mathfrak{p}^{-1}\xi)E(\xi)$, where $E(\xi)$ is the unitary matrix defined in {\rm (\ref{eqn:matE})}.
\end{proposition}
\proof The $(r,s)$-th block of the matrix $P(\mathfrak{p}^{-1}\xi)E(\xi)$ is the matrix
\[
\sum\limits_{t=0}^{q-1} P^{rt}(\mathfrak{p}^{-1}\xi)E^{ts}(\xi).
\]
The $(l,j)$-th entry in this block is equal to
\begin{eqnarray*}
&  & \sum\limits_{t=0}^{q-1}\sum\limits_{m=1}^N p^{rt}_{lm}(\mathfrak{p}^{-1}\xi) {\mathcal E}^{ts}_{mj}(\xi) \\
& = & \sum\limits_{t=0}^{q-1}\sum\limits_{m=1}^N\sum\limits_{k\in\N_0} h^r_{l,m,qk+t} \overline{\chi_k(\mathfrak{p}^{-1}\xi)}\delta_{mj}q^{-1/2}\overline{\chi\bigl(u(t)(\xi+\mathfrak{p}u(s))
\bigr)} \\
& = & \sum\limits_{t=0}^{q-1}\sum\limits_{k\in\N_0} h^r_{l,m,qk+t}\overline{\chi_k(\mathfrak{p}^{-1}\xi)}q^{-1/2}\overline{\chi\bigl(u(t)(\xi+\mathfrak{p}u(s))\bigr)}.
\end{eqnarray*}
Now, the $(l,j)$-th entry in the $(r,s)$-th block of $H(\xi)$ is
\begin{eqnarray*}
& & h^r_{lj}(\xi+pu(s))\\
& = & q^{-1/2}\sum\limits_{k\in\N_0} h^r_{ljk}\overline{\chi\bigl(u(k)(\xi+\mathfrak{p}u(s))\bigr)} \\
& = & q^{-1/2}\sum\limits_{t=0}^{q-1}\sum\limits_{k\in\N_0} h^r_{l,j,qk+t}\overline{\chi\bigl(u(qk+t)(\xi+\mathfrak{p}u(s))\bigr)} \\
& = & q^{-1/2}\sum\limits_{t=0}^{q-1}\sum\limits_{k\in\N_0} h^r_{l,m,qk+t}\overline{\chi(\mathfrak{p}^{-1}u(k)\xi+u(k)u(s)+u(t)\xi+\mathfrak{p}u(t)u(s))} \\
& = & q^{-1/2}\sum\limits_{t=0}^{q-1}\sum\limits_{k\in\N_0} h^r_{l,m,qk+t}\overline{\chi_k(\mathfrak{p}^{-1}\xi)}\overline{\chi\bigl(u(t)(\xi+\mathfrak{p}u(s))\bigr)}.
\end{eqnarray*}
\qed

In particular, we have
\begin{equation*}
H^*(\xi)H(\xi)=E^*(\xi)P^*(\mathfrak{p}^{-1}\xi)P(\mathfrak{p}^{-1}\xi)E(\xi).
\end{equation*}
Since $E(\xi)$ is unitary by Lemma \ref{lem:matrix},
$H^*(\xi)H(\xi)$ and $P^*(\mathfrak{p}^{-1}\xi)P(\mathfrak{p}^{-1}\xi)$ are similar matrices.

Let $\lambda(\xi)$ and $\Lambda(\xi)$ respectively be the maximal and minimal eigenvalues of the positive definite matrix $H^*(\xi)H(\xi)$, and let $\lambda=\inf\limits_\xi\lambda(\xi)$ and $\Lambda=\sup\limits_\xi \Lambda (\xi)$.  Suppose $0<\lambda\leq \Lambda<\infty$. Then we have
\[
\lambda I \leq H^*(\xi)H(\xi) \leq \Lambda I\quad{\rm for~a.e.}~ \xi\in\mathfrak{D}.
\]
This is equivalent to say that
\[
\lambda I \leq P^*(\xi)P(\xi) \leq \Lambda I\quad{\rm for~a.e.}~ \xi\in\mathfrak{D}.
\]
Then by Lemma \ref{lem:frame}, for all $g\in L^2(K)$, we have
\begin{eqnarray}\label{eqn:lambda1}
\lambda\sum\limits_{s=0}^{q-1}\sum\limits_{l=1}^N\sum\limits_{k\in\N_0} \left|\ip{g}{\varphi^{(s)}_l(\cdot-u(k))}\right|^2 & \leq &
\sum\limits_{s=0}^{q-1}\sum\limits_{l=1}^N\sum\limits_{k\in\N_0} \left|\ip{g}{\psi^s_l(\cdot-u(k))}\right|^2 \nonumber \\
& \leq & \Lambda\sum\limits_{s=0}^{q-1}\sum\limits_{l=1}^N\sum\limits_{k\in\N_0} \left|\ip{g}{\varphi^{(s)}_l(\cdot-u(k))}\right|^2,
\end{eqnarray}
where $\varphi^{(s)}_l$ is defined in (\ref{eqn:phisj}).
Since
\begin{equation*}
\sum\limits_{l=1}^N\sum\limits_{k\in\N_0}\left|\ip{g}{q^{1/2}\varphi_l(\mathfrak{p}^{-1}\cdot-u(k))}\right|^2
=\sum\limits_{s=0}^{q-1}\sum\limits_{l=1}^N\sum\limits_{k\in\N_0}\left|\ip{g}{\varphi^{(s)}_l(\cdot-u(k))}\right|^2,
\end{equation*}
which follows from (\ref{eqn:phisj}), inequality (\ref{eqn:lambda1}) can be written as
\begin{eqnarray}\label{eqn:lambda2}
\lambda \sum\limits_{l=1}^N\sum\limits_{k\in\N_0}\left|\ip{g}{q^{1/2}\varphi_l(\mathfrak{p}^{-1}\cdot-u(k))}\right|^2
& \leq & \sum\limits_{s=0}^{q-1}\sum\limits_{l=1}^N\sum\limits_{k\in\N_0} \left|\ip{g}{\psi^s_l(\cdot-u(k))}\right|^2 \nonumber \\
& \leq & \Lambda\sum\limits_{l=1}^N\sum\limits_{k\in\N_0}\left|\ip{g}{q^{1/2}\varphi_l(\mathfrak{p}^{-1}\cdot-u(k))}\right|^2.
\end{eqnarray}

This is the \emph{splitting trick} for frames.

We now apply the splitting trick to the functions $\{\psi^s_l:1\leq l\leq N\}$ for each $s$, $0\leq s\leq q-1$. We have
\begin{eqnarray}
\lambda\sum\limits_{l=1}^N\sum\limits_{k\in\N_0} \left|\ip{g}{q^{1/2}\psi^s_l(\mathfrak{p}^{-1}\cdot-u(k))}\right|^2
& \leq & \sum\limits_{r=0}^{q-1}\sum\limits_{l=1}^N\sum\limits_{k\in\N_0} \left|\ip{g}{\psi^{s,r}_l(\cdot-u(k))}\right|^2  \nonumber \\
& \leq & \Lambda\sum\limits_{l=1}^N\sum\limits_{k\in\N_0} \left|\ip{g}{q^{1/2}\psi^s_l(\mathfrak{p}^{-1}\cdot-u(k))}\right|^2,\label{eqn:split1}
\end{eqnarray}
where $\psi^{s,r}_l,0\leq r\leq q-1$ are defined as:
\begin{equation}\label{eqn:fsrl.frame}
\psi^{s,r}_l(x)= \sum\limits_{l=1}^N\sum\limits_{k\in\N_0}h^s_{ljk}q^{1/2}\psi^r_j(\mathfrak{p}^{-1}x-u(k));
~0\leq s\leq q-1,1\leq l\leq N.
\end{equation}
Summing (\ref{eqn:split1}) over $0\leq s\leq q-1$, we have
\begin{eqnarray*}
\lambda\sum\limits_{s=0}^{q-1}\sum\limits_{l=1}^N\sum\limits_{k\in\N_0} \left|\ip{g}{q^{1/2}\psi^s_l(\mathfrak{p}^{-1}\cdot-u(k))}\right|^2
& \leq & \sum\limits_{s=0}^{q-1}\sum\limits_{r=0}^{q-1}\sum\limits_{l=1}^N\sum\limits_{k\in\N_0} \left|\ip{g}{\psi^{s,r}_l(\cdot-u(k))}\right|^2 \\
& \leq & \Lambda\sum\limits_{s=0}^{q-1}\sum\limits_{l=1}^N\sum\limits_{k\in\N_0} \left|\ip{g}{q^{1/2}\psi^s_l(\mathfrak{p}^{-1}\cdot-u(k))}\right|^2 .
\end{eqnarray*}
Using (\ref{eqn:lambda2}), we obtain
\begin{eqnarray}\label{eqn:split2}
\lambda^2 \sum\limits_{l=1}^N\sum\limits_{k\in\N_0}\left|\ip{g}{q^{2/2}\varphi_l(\mathfrak{p}^2\cdot-u(k))}\right|^2
& \leq & \sum\limits_{s=0}^{q-1}\sum\limits_{r=0}^{q-1}\sum\limits_{l=1}^N\sum\limits_{k\in\N_0} \left|\ip{g}{\psi^{s,r}_l(\cdot-u(k))}\right|^2 \nonumber \\
& \leq & \Lambda^2 \sum\limits_{l=1}^N\sum\limits_{k\in\N_0}\left|\ip{g}{q^{2/2}\varphi_l(\mathfrak{p}^2\cdot-u(k))}\right|^2.
\end{eqnarray}

We now define the wavelet frame packets similar to the orthonormal case. We start with the functions $\varphi_1, \varphi_2, \dots, \varphi_N$. Apply the splitting trick to the space
\[
\overline{{\rm span}}\{q^{1/2}\varphi_l(\mathfrak{p}^{-1}\cdot-u(k)):1\leq l\leq N, k\in\N_0\}
\]
to get the functions $\{\psi_l^s: 1\leq l\leq N, 0\leq s\leq q-1\}$ (see~(\ref{eqn:lambda2})). Now for any integer $n\geq 0$, we define $\psi^n_l$, $1\leq l\leq N$, recursively as follows. Suppose that $\psi^r_l$ is already defined for $r\in\N_0$ and $1\leq l\leq N$. Then for $0\leq s\leq q-1$ and $1\leq l\leq N$, define
\[
\psi_l^{s+qr} = \sum\limits_{j=1}^N\sum\limits_{k\in\N_0}h^s_{ljk}q^{1/2}\psi^r_j(\mathfrak{p}^{-1}\cdot-u(k)).
\]

Comparing this with equation~(\ref{eqn:fsrl.frame}), we see that
\begin{eqnarray*}
\{\psi^{s,r}_l:0\leq r,s\leq q-1\} & = & \{\psi^{s+qr}_l:0\leq r,s\leq q-1\} \\
& = & \{\psi^n_l:0\leq n\leq q^2-1\}.
\end{eqnarray*}
So (\ref{eqn:split2}) can be written as
\begin{eqnarray*}
\lambda^2\sum\limits_{l=1}^N\sum\limits_{k\in\N_0}\left|\ip{g}{q^{2/2}\varphi_l(\mathfrak{p}^{-2}\cdot-u(k))}\right|^2
& \leq & \sum\limits_{n=0}^{q^2-1}\sum\limits_{l=1}^N\sum\limits_{k\in\N_0} \left|\ip{g}{\psi^n_l(\cdot-u(k))}\right|^2 \nonumber \\
& \leq & \Lambda^2 \sum\limits_{l=1}^N\sum\limits_{k\in\N_0}\left|\ip{g}{q^{2/2}\varphi_l(\mathfrak{p}^{-2}\cdot-u(k))}\right|^2.
\end{eqnarray*}
By induction, we get for each $j\geq 1$
\begin{eqnarray}\label{eqn:splitj}
\lambda^j\sum\limits_{l=1}^N\sum\limits_{k\in\N_0}\left|\ip{g}{q^{j/2}\varphi_l(\mathfrak{p}^{-j}\cdot-u(k))}\right|^2
& \leq & \sum\limits_{n=0}^{q^j-1}\sum\limits_{l=1}^N\sum\limits_{k\in\N_0} \left|\ip{g}{\psi^n_l(\cdot-u(k))}\right|^2 \nonumber \\
& \leq & \Lambda^j\sum\limits_{l=1}^N\sum\limits_{k\in\N_0}\left|\ip{g}{q^{j/2}\varphi_l(\mathfrak{p}^{-j}\cdot-u(k))}\right|^2.
\end{eqnarray}

We summarize the above discussion in the following theorem.
\begin{theorem}
Let $\{\varphi_l:1\leq l\leq N\}\subset L^2(K)$ be such that $\{\varphi_l(\cdot-u(k)):1\leq l\leq N,k\in\N_0\}$ is a frame for its closed linear span $V_0$, with frame bounds $C_1$ and $C_2$ . Let $H(\xi),H_r(\xi),\lambda$ and $\Lambda $ be as above. Assume that all entries of $H(\xi)$ are bounded measurable functions such that $0< \lambda\leq \Lambda <\infty$. Let $\{\psi^n_l:n\geq 0,1\leq l\leq N\}$ be the wavelet frame packets and let
$V_j=\{f\in L^2(K): f(\mathfrak{p}^{j}\cdot)\in V_0\}$. Then for all $j\geq 0$, the system of functions
\[
\{\psi^n_l(\cdot-u(k)): 0\leq n\leq q^j-1, 1\leq l\leq N, k\in\N_0\}
\]
is a frame of $V_j$ with frame bounds $\lambda^jC_1$ and  $\Lambda^jC_2$.
\end{theorem}
\proof Since $\{\varphi_l(\cdot-u(k)):1\leq l\leq N.k\in\N_0\}$ is a frame of $V_0$ with frame bounds $C_1$ and $C_2$, it is clear that for all $j$
\[
\{q^{j/2}\varphi_l(\mathfrak{p}^{-j}\cdot-u(k)): 1\leq l\leq N, k\in\N_0\}
\]
is a frame of $V_j$ with the same bounds. So from (\ref{eqn:splitj}), we have
\begin{equation*}
\lambda^jC_1\|g\|^2 \leq \sum\limits_{n=0}^{q^j-1}\sum\limits_{l=1}^N\sum\limits_{k\in\N_0} \left|\ip{g}{\psi^n_l(\cdot-u(k))}\right|^2 \leq \Lambda^jC_2\|g\|^2
\end{equation*}
for all $g\in V_j$.
\qed


\end{document}